\documentstyle[11pt,leqno]{article}
\newtheorem{prop}{Proposition}[section]
\newtheorem{th}[prop]{Theorem}
\newtheorem{cor}[prop]{Corollary}
\newtheorem{lm}[prop]{Lemma}

\newcommand{\bsquare}{\hbox{\rule{6pt}{6pt}}}
\newcommand{\proof}[1]{\noindent{\bf Proof}\hspace{0.3cm}{#1}\hfill\bsquare
\vspace{0.5cm}\par}

\newcommand{\Ker}{{\rm{Ker}}~}
\renewcommand{\Im}{{\rm{Im}}~}

\renewcommand{\H}{{\rm{H}}}

\newcommand{\rank}{{\rm{rank}}~}
\newcommand{\Pic}{{\rm{Pic}}}

\renewcommand{\O}{{\cal{O}}}
\newcommand{\U}{{\cal{U}}}
\newcommand{\G}{{\bf{G}}}

\def\spmapright#1{\smash{%
   \mathop{\hbox to 1.3cm{\rightarrowfill}}
       \limits^{#1}}}
\def\sbmapright#1{\smash{%
   \mathop{\hbox to 1.3cm{\rightarrowfill}}
       \limits_{#1}}}

\def\rmapdown#1{\Big\downarrow
   \rlap{$\vcenter{\hbox{$\scriptstyle#1$}}$ }}

\newcommand{\mapright}[1]{%
\smash{\mathop{%
   \hbox to 1cm{\rightarrowfill}}\limits^{#1}}}
\newcommand{\mapleft}[1]{%
\smash{\mathop{%
   \hbox to 1cm{\leftarrowfill}}\limits^{#1}}}

\begin{document}
  \pagenumbering{arabic}

\title{Formal Brauer groups and the moduli\\
 of abelian surfaces}

\author{G. van der Geer and T. Katsura}

\maketitle

\section{Introduction}
Let $X$ be an algebraic surface over an algebraically closed field
$k$ of characteristic $p > 0$.  We denote by $\Phi_X$ the formal Brauer group
of $X$ and  by $h = h(\Phi_X)$ the height of $\Phi_X$. In a previous paper,
\cite{GK}, we examined the structure of the stratification given by the
height $h$ in the moduli space of K3 surfaces, and we determined the
cohomology class of each stratum. In this paper, we apply the methods of
\cite{GK} to treat the case of abelian surfaces. In this case, the situation
is more  concrete, and so we can more easily determine the structure of  the
stratification given by the height $h(\Phi_{A})$ in the moduli of  abelian
surfaces. For the local structure we refer to \cite{Sh}. 

On the moduli of
principally polarized abelian varieties in positive characteristic there is
another natural stratification, called the Ekedahl-Oort stratification, cf.\
\cite{O1}, and one can calculate the corresponding cycle classes
\cite{G2}.  Although our three strata coincide set-theoretically with strata
of the Ekedahl-Oort stratification, there is a subtle difference: one of the
strata comes with multiplicity $2$.  

We will here summarize our results. We consider the moduli stack $M={\cal
A}_2$ of principally polarized abelian surfaces over $k$; alternatively,
we can consider the moduli spaces $M={\cal A}_{2, n}$
$(n \geq 3, ~ p \not| n)$   of principally polarized abelian
surfaces  with level $n$-structure. We know that $M$ is a 3-dimensional
algebraic stack (variety). We let  $\pi:{\cal X} \rightarrow M$ be
the universal family  over $M$. 
We set
$$
M^{(h)}:=\{ s \in M: h(\Phi_{{\cal X}_s}) \geq h\}.
$$
Note that $M^{(3)} = M^{(\infty)}$. The moduli stack $M$ possesses a natural
compactification $\tilde{M}$ which is an example of a smooth toroidal
compactification, cf.\ \cite{FC} or \cite{Al}. It carries a universal family 
${\tilde{\pi}}: \tilde{{\cal X}} \rightarrow \tilde{M}$. We can extend in a
natural way the loci $M^{(h)}$ to loci  in $\tilde{M}$ which are again
denoted by ${M}^{(h)}$.

For $h = 2, 3$, we denote by $M^{(h)}_{F}$ the scheme-theoretic zero
locus in $\tilde{M}$ of the Frobenius action
$$
F : \H^{2}(A, W_{h - 1}({\O}_{A})) \longrightarrow
\H^{2}(A, W_{h - 1}({\O}_{A})).
$$
(See Section 7 for details.) We set
$M^{(\infty)}_{F}= M^{(3)}_{F}$.
We will show $M^{(2)} = M^{(2)}_{F}$ and
$M^{(\infty)} =( M^{(\infty)}_{F})_{red}$ in $\tilde{M}$, and the following
theorem  describes the Chow classes of these loci in $\tilde{M}$.
We denote
by $v$ the first Chern class of the coherent sheaf
$R^0{\tilde{\pi}}\Omega^2_{\tilde{{\cal X}}/\tilde{M}}$ in the Chow group
$CH_{\bf Q}^{1}(\tilde{M})$.

{\bf Theorem.}
{\sl The classes of the loci} ${M}_F^{(h)}$ {\sl in the Chow group} $CH_{\bf
Q}^{*}({\tilde{M}})$
{\sl are  given by}
$$
     M^{(2)}_{F} = (p-1)v, ~M^{(\infty)}_{F} = (p-1)(p^2-1)v^{2}.
$$

To prove these results, we use the characterization of the height of
$\Phi_A$ by the action of the Frobenius morphism on $\H^2(W_i(\O_A))$
which was obtained in \cite{GK} (cf. Section 2). We also investigate
the natural images of $\H^{1}(A, B_{i})$
and $\H^{1}(A, Z_{i})$ in $\H^{1}(A, \Omega_{A}^{1})$, which are related
to the tangent space to $M^{(h)}$. (For the case of K3 surface, see \cite{Og2}.)

Comparing our results with those in \cite{G2}, we have the following
relation  between ${M}^{(h)}$ and $M^{(h)}_{F}$.

{\bf Theorem.}
{\sl In the Chow ring} $CH_{\bf Q}^{*}({\tilde{M}})$ {\sl we have}
$$
M^{(2)}_{F} = M^{(2)}, ~M^{(\infty)}_{F} = 2M^{(\infty)}.
$$

To examine the structure of abelian surfaces $A$, one usually employs the
first cohomology group $\H^{1}(A, \O_{A})$. Instead, in this paper, we
mainly use $\H^{2}(A, W_{i}(\O_{A}))$ and we will show that our techniques
developed in  \cite{GK}  also work for the study of abelian surfaces.      

Most part of this work was done in Max-Planck-Institut f\"ur Mathematik 
in Bonn. The authors would like to thank Max-Planck-Institut for
excellent working conditions during their stay in 1998/99.
\section{Preliminaries}
Let $X$ be a non-singular complete algebraic variety defined over an
algebraically closed field $k$ of characteristic $p>0$. We denote by
$W_n(\O_X)$ the sheaf of Witt vectors of length $n$ 
(cf. J.-P.\ Serre, \cite{S}). The sheaf $W_n(\O_X)$  is a coherent 
sheaf of rings. Sometimes, we write $W_{n}$ instead of $W_n(\O_X)$
for the sake of simplicity.  We denote by $F$ (resp. $V$, resp. $R$) the 
Frobenius map (resp. the Verschiebung, resp. the restriction map). They
satisfy relations 
$$
RVF=FRV=RFV=p.
$$
The cohomology groups ${\H}^i(X,W_n(\O_X))$ are finitely generated 
$W_n(k)$-modu\-les. The projective system $\{ W_n(\O_X), R\}_{n=1,2,\ldots}$
defines the cohomology group
$$
\H^i(X,W(\O_X)) = \hbox{\rm proj.\ lim } \H^i(X,W_n(\O_X)).
$$
This is a $W(k)$-module, but not necessarily a finitely
generated $W(k)$-module. The semi-linear operators $F$
and $V$ act on it with relations $FV=VF=p$.

For any   artinian local $k$-scheme $S$ with residue field $k$,
we consider the kernel
$$
\Phi^i_{X}(S)=\Ker \H^i(X\times S, \G_m) \longrightarrow H^i(X,\G_m).
$$
Here $\G_m$ is the multiplicative group of dimension 1 and
cohomology is \'etale cohomology. This gives
a contravariant functor $\Phi^i_{X} : {\rm Art} \to {\rm Ab}$
from the category of artinian local $k$-schemes with residue field $k$
to the category of abelian groups. This functor was introduced
and investigated  by Artin and Mazur \cite{AM}. They showed that the tangent
space of
$\Phi^{i}_{X}$ is given by
$$
T_{\Phi^{i}_{X}}= \H^{i}(X,{\O}_{X}),
$$
and proved the following crucial theorem.
\begin{th}\label{th;AM}
The Dieudonn\'e module of the Artin-Mazur formal
group $\Phi^{i}_{X}$ is given by
$$
D(\Phi^{i}_{X}) \cong \H^{i}(X,W({\O}_X)).
$$
\end{th}

If $A$ is an abelian surface over $k$ then $\Phi^{1}_{A}$ is
isomorphic to the formal completion of the Picard scheme $\Pic^{0} (A)$ 
of $A$, which is a formal Lie group. Moreover, we have $\H^{3}(A, {\O}_{A}) = 0$
and $\H^{2}(A, {\O}_{A}) \cong k$ .
Therefore, by a criterion of pro-representability in \cite{AM}, we know
that $\Phi_{A} = \Phi^{2}_{A}$ is pro-representable by a formal Lie group
of dimension 1, which is called a formal Brauer group of $A$.
One-dimensional formal groups are classified by their height $h$, which is  a
positive integer or $h=\infty$ for the case of the additive group.
\section{Abelian Surfaces}
For an abelian surface $A$, we denote by $[p]_{A}$ the homomorphism
given by multiplication by $p$ on $A$, and by $\Ker[p]_{A,
red}$ the reduced group scheme associated with the kernel $\Ker[p]_{A}$. 
The dimension of $\Ker[p]_{A, red}$ over ${\bf F}_{p}$ is called the $p$-rank
of $A$, and we denote it by $r(A)$. As is well known, we have $0 \leq r(A)
\leq 2$. Since any regular 1-form on an abelian surface is d-closed, the
Cartier  operator  $C$ acts on $\H^{0}(A, \Omega_{A}^{1})$. We have also a
natural action of  the Frobenius map on $\H^{1}(A, \O_{A})$, which is dual
to the Cartier  operator $C$ on $\H^{0}(A, \Omega_{A}^{1})$. The abelian
surface
$A$ is said to be {\sl ordinary} if $r(A) = 2$, which is equivalent to 
 $F$ (resp. $C$) being bijective on $\H^{1}(A, \O_{A})$ 
(resp. $\H^{0}(A, \Omega_{A}^{1})$). Furthermore, $A$ is said to be
{\sl supersingular} if $A$ is  isogenous to a product of two supersingular
elliptic curves, which is  equivalent to the effect that $F$ (resp. $C$) is
nilpotent on $\H^{1}(A, \O_{A})$  (resp. $\H^{0}(A, \Omega_{A}^{1})$).
Finally,  $A$ is said to be {\sl superspecial} if $A$ is  isomorphic to a
product of two supersingular elliptic curves, which is  equivalent to 
 $F$ (resp. $C$) being the zero map on $\H^{1}(A, \O_{A})$  (resp.
$\H^{0}(A,
\Omega_{A}^{1})$). The superspecial case can be characterized numerically
by the $a$-number of $A$. Here $a(A)= \dim_k {\rm Hom}(\alpha_p,A)$ and we
know 
\begin{itemize}
\item{} $a=0 \iff r=2$, 
\item{} $a\geq 1\iff r\leq 1$,
\item{} $a=2 \iff A$ is superspecial.
\end{itemize}
For details on abelian surfaces,  see
\cite{Mu1},
\cite{O1} and \cite{OO}. The following lemma is well-known.

\begin{lm}\label{lm:abelian}
The height $h$ of the formal Brauer group 
$\Phi_{A}$ of an abelian surface $A$  is  as follows:
\begin{itemize}
\item[$({\rm 1})$] $h = 1$ if $r(A)=2$, i.e. $A$ is ordinary,
\item[$({\rm 2})$] $h = 2$ if $r(A) = 1$,
\item[$({\rm 3})$] $h = \infty$ if $r=0$, i.e. $A$ is supersingular.
\end{itemize}
\end{lm}
\proof{We denote by $\H^{i}_{cris}(A)$ the i-th cristalline cohomology of
$A$ and as usual by $\H^{i}_{cris}(A)_{[\ell, \ell  + 1[}$ the
additive group of elements in $\H^{i}_{cris}(A)$ whose slopes are in the
interval $[\ell,\ell + 1[$. By the general theory
in Illusie \cite{I}, we have 
$$
 \quad \H^{2}(A,W({\O}_A)) \otimes_{W}K \cong 
(\H^{2}_{cris}(A) \otimes_{W}K)_{[0, 1[}
$$
with $K$ the quotient field of $W$. The theory of Dieudonn\'e modules
implies 
$$
 h = \dim_{K} D(\Phi_{A}) = \dim_{K} \H^{2}(A,W({\O}_A)) \otimes_{W}K \quad 
\mbox{if} ~h < \infty,
$$
and $\dim_{K} D(\Phi_{A}) = 0$ if $h = \infty$.
We know the slopes of $\H^{1}_{cris}(A)$ for each case. Since
we have 
$$
  \H^{2}_{cris}(A) \cong \wedge^{2} \H^{1}_{cris}(A),
$$
counting the number of slopes in $[0, 1[$ of 
$\H^{2}_{cris}(A)$ gives the result.}
 
We shall
need the following lemma. 

\begin{lm}\label{lm:exact}
For an abelian surface $A$ the following sequence is exact.
$$
0 \rightarrow \H^{2}(A, W_{n-1}(\O_{A})) \stackrel{V}{\longrightarrow}  
\H^{2}(A, W_{n}(\O_{A}))
\stackrel{R^{n-1}}{\longrightarrow} \H^{2}(A, \O_{A}) \rightarrow  0.
$$
\end{lm}
\proof{By the exact sequence
$$
0 \rightarrow W_{n-1}(\O_{A})\stackrel{V}{\longrightarrow}  
W_n(\O_A) \stackrel{R^{n-1}}{\longrightarrow} \O_A \rightarrow 0,
$$  
we have the long exact sequence
$$
 \rightarrow \H^{2}(A, W_{n-1}(\O_{A})) \stackrel{V}{\longrightarrow} 
\H^{2}(A, W_{n}(\O_{A}))
\stackrel{R^{n-1}}{\longrightarrow} \H^{2}(A, \O_{A}) \rightarrow 0.
$$
Since the Picard scheme of $A$ is reduced, all Bockstein operaters vanish
by Mumford \cite{Mu2}, p.\ 196. The result follows from this fact.}

Using this lemma and $FV = VF$, we have the following lemma.

\begin{lm}\label{lm:ker}
If $F$ acts as  zero on $\H^{2}(A, W_{n}(\O_{A}))$ then 
$F$ acts as  zero on $\H^{2}(A, W_{i}(\O_{A}))$ for any $i < n$.
\end{lm}

The following theorem is an analogue of a result in \cite{GK}.
\begin{th}\label{th:character}
The height satisfies  $h(\Phi_A) \geq i+1$  if and
only if the Frobenius map $ F :  H^2(A,W_i(\O_A)) \rightarrow
H^2(A,W_i(\O_A))$ is the zero map. In particular, we have the following
characterization of the height:
$$
h(\Phi_A) = \min \{ i \geq 1~ \mid~ [F: \H^2(W_i(\O_A)) \rightarrow
\H^2(W_i(\O_A))] \neq 0 \}.
$$
\end{th}

\begin{cor} Let $A$ be an abelian surface over $k$. Then 
$A$ is supersingular if and only if  the Frobenius endomorphism  
$F: \H^2(A,W_{2}(\O_A)) \rightarrow \H^2(A,W_{2}(\O_A))$ is zero.
\end{cor}
\proof{If the height is finite, then we know $h\leq 2$ by 
Lemma \ref{lm:abelian}. So this corollary follows from Theorem
\ref{th:character} and  Lemma \ref{lm:abelian} .}

\begin{cor}
Let $A$ be an abelian surface over $k$.
Set $\H= \H^2(A,W_{2}(\O_A))$. If the height $h$ of $\Phi_A$ is finite, then
$F(\H)= R^{h-1}V^{h-1}(\H)$.
\end{cor}
\proof{The $(h-1)$-th step $V^{h-1}\H^2(W(\O_A))$ in the filtration
$$
V^{2}\H^2(W(\O_A)) \subset V\H^2(W(\O_A)) \subset \H^2(W(\O_A))
$$
maps surjectively to the corresponding step $R^{h-1}V^{h-1}\H$ of the
filtration on $\H$. Therefore, this corollary follows from
$$
V^{h-1}\H^2(W(\O_A)) = F\H^2(W(\O_A)).
$$ 
}

For an element $\omega \in \H^2(W(\O_A))$, we denote by $\bar\omega$
the natural restriction of $\omega$ in $\H^2(A,W_i(\O_A))$. The following
corollary follows immediately from $FV=VF$.

\begin{cor}\label{cor:fundamental}
Let $A$ be an abelian surface over $k$.
If $h(\Phi_A)=h < \infty$ and if
$\{ \omega, V^{h-1}\omega\} $ is a $W$-basis of
$\H^2(X,W(\O_X))$ then $F$ acts as zero on $\H^2(A,W_i(\O_A))$ if and only if
$F(\bar\omega)=0$, with $\bar\omega$ the image of $\omega$ in
$\H^2(A,W_i(\O_A))$.
\end{cor}

\begin{cor}\label{cor: KerF}
Let $A$ be an abelian surface over $k$. Putting $h(\Phi_A)=h$,
we have
$$
\dim_{k} \Ker[F: \H^2(W_i) \rightarrow \H^2(W_i)] =
\left\{
\begin{array}{l}
 0 \quad \mbox{if}~h = 1 \\
  1 \quad \mbox{if}~ h = 2 \\
 i  \quad \mbox{if}~h = \infty
\end{array}
\right.
$$
\end{cor}
\proof{If $h = 1$ then $F : \H^{i}(A, \O_{A}) \rightarrow \H^{i}(A, \O_{A})$
is a $p$-linear isomorphism for any $i$. Therefore, using exact sequences
in Lemma \ref{lm:exact} we can inductively show that $F: \H^2(W_i)
\rightarrow \H^2(W_i)$ is an isomorphism for any $i$.

If $h = 2$, then
by Corollary \ref{cor:fundamental} we have 
$\dim_{k} \Ker[F: \H^2(\O_{A}) \rightarrow \H^2(\O_{A})] = 1$.
Assume $i \geq 2$. Using the notation in 
Corollary \ref{cor:fundamental}, we know that 
$\langle V^{i-1}R^{i-1}{\bar\omega} \rangle$ 
is a basis of $\Ker F$. Therefore, it is one-dimensional.

If $h = \infty$, then by Theorem \ref{th:character}, we see that 
$\Ker[F: \H^2(W_i) \rightarrow \H^2(W_i)] = \H^2(W_i)$. Therefore, using 
Lemma \ref{lm:exact} we get the result by induction.}

Let $A$ be an abelian surface  and assume that $F$ is zero on $\H^2(\O_A)$.
Then we have $F\H^2(W_2) \subset V\H^2(\O_A)$ and
$F$ vanishes on $V\H^{2}(\O_A)$. Since we have the natural 
$({\sigma}^{-1})$-isomorphism $\H^2(\O_A) \cong V\H^2(\O_A)$, we have the
induced homomorphism 
$$\phi_{2} : \H^{2}(\O_{A}) \cong \H^2(W_2)/V\H^2(\O_{A})
\rightarrow V\H^2(\O_A) \cong  \H^2(\O_ A).
$$
This map is ${\sigma}^{2}$-linear.
The following theorem is clear by the construction.
\begin{th}\label{th:h=2}
 Suppose $F$ is zero on $\H^2(\O_{A})$. Then $F$ is zero on $\H^2(W_2)$ 
if and only if $\phi_{2} = 0$.
\end{th}

\section{Cohomology Groups of Abelian Surfaces }
Let $A$ be an abelian surface defined over $k$.
We define sheaves $B_i\Omega_A^{1}$ inductively by $B_0\Omega_A^1=0$,
$B_1=dO_A$ and $C^{-1}(B_i\Omega_A^1) = B_{i+1}\Omega_A^1$. 
Similarly, we define sheaves $Z_i\Omega_A^1$ inductively by
$Z_0\Omega_A^1=\Omega_A^1$, $Z_1\Omega_A^1= \Omega^1_{A,{\rm closed}}$, the
sheaf of $d$-closed forms and by setting 
$$
Z_{i+1}\Omega_A^1 : = C^{-1}(Z_i\Omega_A^1).
$$
Sometimes we simply write $B_i$ (resp. $Z_i$) instead of $B_i\Omega_A^{1}$
(resp. $Z_i\Omega_A^1$).
Note that we have the inclusions
$$
0=B_0 \subset B_1 \subset \ldots \subset B_i \subset \ldots \subset Z_i
\subset \ldots \subset Z_{1} \subset Z_0 = \Omega_X^1.
$$
The sheaves $B_i$ and $Z_i$ can be viewed as locally free subsheaves of
$(F^i)_*\Omega_A^1$ on $A^{(p^i)}$.  
The inverse Cartier operator gives rise to
an isomorphism 
$$ 
C^{-i} : \Omega_A^1 \stackrel{\sim}{\longrightarrow}
Z_i\Omega_A^1/B_i\Omega_A^1.
$$
We also have an exact sequence   
$$ 
0 \rightarrow Z_{i+1}\Omega_A^1 \longrightarrow Z_i\Omega_A^1 
\stackrel{dC^i}{\longrightarrow} d\Omega_A^1 \rightarrow 0. 
$$
\begin{lm}\label{lm:Bi-0}
Let  $A$ be an abelian surface over $k$. Then the dimension $h^j(A,B_i)$
for $i\geq 1$ is as in the following table.
\end{lm} 

\font\tablefont=cmr8
\def\quad{\hskip 0.6em\relax}
\def\quod{\hskip 0.6em\relax}
\def\vhop{
    height2pt&\omit&&\omit&&\omit&&\omit&\cr}
$$
\vcenter{
\tablefont
\lineskip=1pt
\baselineskip=10pt
\lineskiplimit=0pt
\setbox\strutbox=\hbox{\vrule height .7\baselineskip
                                depth .3\baselineskip width0pt}%
\offinterlineskip
\hrule
\halign{&\vrule#&\strut\quod\hfil#\quad\cr
\vhop
&${\rm type}$&&$h^0(B_i)$&&$h^1(B_i)$&&$h^2(B_i)$&\cr
\vhop
\noalign{\hrule}
\vhop
&$h=1$&&$0$&&$0$&&$0$&\cr
\vhop
\noalign{\hrule}
\vhop
&$h=2$&&$1$&&$2$&&$1$&\cr
\vhop
\noalign{\hrule}
\vhop
&$h=\infty$, $a=1$&&${\cases{1 & $i=1$\cr 2 & $i\geq
2$}}$&&${
\cases{2 & $i=1$\cr 2+i & $i\geq 2$}}$&&$i$&\cr
\vhop
\noalign{\hrule}
\vhop
&$h=\infty$, $a=2$&&$2$&&$2+i$&&$i$&\cr
\vhop
\noalign{\hrule}
\vhop
}
\hrule
}
$$
\proof{The $h^0$-part of this lemma follows from  the 
natural inclusion 
$$
    \H^{0}(A, B_{i}) \rightarrow \H^{0}(A, \Omega^{1}_{A})
$$
and by considering the action of Cartier operator in each case.
For the $h^1$-part we now consider the map $D_i: W_i(\O_A) \rightarrow
\Omega_A^1$ of  sheaves which was introduced by Serre in the following way: 
$$
D_i(a_0,a_1,\ldots,a_{i-1}) =  a_0^{p^{i-1}-1}da_0 +\ldots + 
a_{i-2}^{p-1}da_{i-2} +
da_{i-1}.
$$
Serre showed that
this map is an additive homomorphism with $D_{i+1}V = D_{i}$, and that
this induces an injective map of sheaves
$$
D_i : W_i(\O_A)/FW_i(\O_A) \longrightarrow \Omega^1_A 
$$
inducing an isomorphism 
$D_i: W_i(\O_A)/FW_i(\O_{A})  \stackrel{\sim}{\longrightarrow}
B_i\Omega_A^1$.
The exact sequence  
$$
0 \rightarrow W_i  \stackrel{F}{\longrightarrow}
W_i \longrightarrow W_i/FW_i \rightarrow 0
$$ 
gives rise to the  exact sequence
$$
\begin{array}{c}
0 \rightarrow \H^{1}(W_i)/F\H^{1}(W_i) \longrightarrow \H^1(W_i/FW_i) \\
\longrightarrow\Ker[F: \H^2(W_i) \longrightarrow \H^2(W_i)] \rightarrow 0.
\end{array}
$$
Now the $W_i(k)$-modules $\H^{1}(W_i)/F\H^{1}(W_i)$ and $\H^{1}(W_i/FW_i)$
are vector spaces over $k \cong W_{i}(k)/pW_{i}(k)$.
We know the dimension of the kernel of $F$ on $\H^2(W_i)$ 
by Corollary \ref{cor: KerF}. 
Since $\H^{1}(W_i)$ is a $W_{i}(k)$-module of finite length  the exact
sequence
$$
\begin{array}{c}
0 \rightarrow \H^{0}(W_i/FW_i) \longrightarrow
\H^{1}(W_i)\stackrel{F}{\longrightarrow} \H^{1}(W_i) \\
 \qquad \longrightarrow \H^{1}(W_i)/F\H^{1}(W_i) 
\rightarrow 0,
\end{array}
$$
implies 
$\dim_{k} \H^{0}(W_i/FW_i) = \dim_{k}\H^{1}(W_i)/F\H^{1}(W_i)$.
Using $\H^{0}(W_i/FW_i) \cong \H^{0}(B_{i})$ and Lemma \ref{lm:Bi-0}
the statement about $h^1$ follows. As to the statement about $h^2$,
consider the exact sequence
$$
0 \rightarrow d{\O_{A}} \longrightarrow B_{i+1} 
\stackrel{C}{\longrightarrow}
B_{i} \rightarrow 0,
$$
and use $\chi (B_{i}) = 0$. Then the result follows immediately.}

The exact sequence
$$
 0 \rightarrow d\Omega_{A}^{1} \longrightarrow \Omega_{A}^{2}
\stackrel{C}{\longrightarrow} \Omega_{A}^{2} \rightarrow 0.
$$ implies:

\begin{lm}\label{lm: dOmega}
Let $A$ be an abelian surface over $k$. Then the dimension
$h^i(d\Omega_A^1)$ of $\H^i(A,d\Omega_A^1)$ is as in the following table.
\end{lm}

\font\tablefont=cmr8
\def\quad{\hskip 0.6em\relax}
\def\quod{\hskip 0.6em\relax}
\def\vhop{
    height2pt&\omit&&\omit&&\omit&&\omit&\cr}
$$
\vcenter{
\tablefont
\lineskip=1pt
\baselineskip=10pt
\lineskiplimit=0pt
\setbox\strutbox=\hbox{\vrule height .7\baselineskip
                                depth .3\baselineskip width0pt}%
\offinterlineskip
\hrule
\halign{&\vrule#&\strut\quod\hfil#\quad\cr
\vhop
&${\rm type}$&&$h^0(d\Omega_A^1)$&&$h^1(d\Omega_A^1)$&&$h^2(d\Omega_A^1)$&\cr
\vhop
\noalign{\hrule}
\vhop
&$h=1$&&$0$&&$0$&&$0$&\cr
\vhop
\noalign{\hrule}
\vhop
&$h=2$&&$1$&&$2$&&$1$&\cr
\vhop
\noalign{\hrule}
\vhop
&$h=\infty$, $a=1$&&$1$&&$2$&&$1$&\cr
\vhop
\noalign{\hrule}
\vhop
&$h=\infty$, $a=2$&&$1$&&$3$&&$2$&\cr
\vhop
\noalign{\hrule}
\vhop
}
\hrule
}
$$

Similarly, the exact sequence
$ 0 \rightarrow B_{i} \longrightarrow Z_{i} 
\stackrel{C^{i}}{\longrightarrow} \Omega_{A}^{1} \rightarrow 0
$
implies
\begin{lm} 
If $h = 1$ we have $h^0(Z_i)=2$, $h^1(Z_i)=4$ and $h^2(Z_i)=2$.
\end{lm}

We also need the following lemma.
\begin{lm}\label{lm:surjective}
Consider the natural inclusions $Z_{i} \hookrightarrow Z_{i-1}$ $(i \geq 1)$
and the induced homomorphisms $\H^{1}(A, Z_{i}) \rightarrow \H^{1}(A,
Z_{i-1})$.
The surjectivity of the map $\H^{1}(A, Z_{i}) \rightarrow \H^{1}(A, Z_{i-1})$
implies the surjectivity of the map $\H^{1}(A, Z_{j}) \rightarrow \H^{1}(A,
Z_{j-1})$  for any $j \geq i$.
\end{lm}
\proof{Suppose that the natural homomorphism 
$\H^{1}(A, Z_{n + 1}) \rightarrow \H^{1}(A, Z_{n})$ is surjective.
By the diagram of exact sequences
$$
\begin{array}{ccccccc}
0 \rightarrow & B_{1} & \longrightarrow & Z_{n + 2} &
\stackrel{C}{\longrightarrow}  & Z_{n + 1} &  \rightarrow 0  \\
      &  \rmapdown{=} &    &\rmapdown{} &    
& \rmapdown{} &                 \\
0 \rightarrow & B_{1} & \longrightarrow & Z_{n + 1} &
\stackrel{C}{\longrightarrow}  & Z_{n} &  \rightarrow 0
\end{array}
$$
we have a diagram of  exact sequences
$$
\begin{array}{cccccccc}
      \rightarrow & \H^{0}(A, B_{1}) & \rightarrow &
\H^{1}(A, Z_{n + 2}) & \stackrel{C}{\longrightarrow} & \H^{1}(A, Z_{n+1})  &
\rightarrow &
\H^{1}(A, B_{1})    \\
  & \rmapdown{=} &  & \rmapdown{} &    &\rmapdown{} &    & \rmapdown{=}   \\
      \rightarrow & \H^{0}(A, B_{1}) & \rightarrow &
\H^{1}(A, Z_{n +1}) & \stackrel{C}{\longrightarrow} & \H^{1}(A, Z_{n}) & 
\rightarrow &
\H^{1}(A, B_{1}).   
\end{array}
$$
>From this diagram, we see that the natural homomorphism 
$\H^{1}(A, Z_{n+2}) \rightarrow \H^{1}(A, Z_{n+1})$
is also surjective. Hence this lemma follows by induction.}
\begin{lm}\label{lm:hikaku}
Assume that $h \geq 2$, and that $A$ is not superspecial. Then,
either  $\dim \H^{1}(A, Z_{i}) = \dim \H^{1}(A, Z_{i-1})$ or 
$\dim \H^{1}(A, Z_{i}) = \dim \H^{1}(A, Z_{i-1}) + 1$.  Moreover, the
latter case occur if and only if the natural homomorphism 
$\H^{1}(A, Z_{i}) \rightarrow \H^{1}(A, Z_{i -1})$ is surjective.
\end{lm}
\proof{By the exact sequence
$
0 \rightarrow  Z_{i}  \longrightarrow  Z_{i-1} 
 \stackrel{dC^{i-1}}{\longrightarrow}  d\Omega_{A}^{1} \rightarrow 0 
$
we have an exact sequence
$$ 
  0 \rightarrow  k  \longrightarrow 
\H^{1}(A, Z_{i})  \longrightarrow  \H^{1}(A, Z_{i-1})  \longrightarrow k. 
$$
The result follows from this exact sequence.}

\begin{prop}
Assume that $h = 2$. Then $\dim \H^{1}(A, Z_{i}) = 4$ for $i \geq 0$.
\end{prop}
\proof{The diagram of exact sequences
$$
\begin{array}{ccccccc}
  & 0  &  & 0 &    & &      \\
  & \downarrow  &  & \downarrow &    & &      \\
  & B_{1}  & =  & B_{1} &    & &      \\
  & \downarrow  &  & \downarrow  &    & &      \\
0 \rightarrow & B_{i}& \longrightarrow &Z_{i} 
&\stackrel{C^{i}}{\longrightarrow} & \Omega_{A}^{1}& \rightarrow 0  \\
      & \rmapdown{C} &    &\rmapdown{C}  &    & \rmapdown{=} &           \\
0 \rightarrow & B_{i-1} & \longrightarrow & Z_{i-1} 
& \stackrel{C^{i-1}}{\longrightarrow} & \Omega_{A}^{1}& \rightarrow 0. \\
  & \downarrow  &  & \downarrow &    & &      \\
  & 0  &  & 0 &    &        &      
\end{array}
$$
gives rise to a diagram of exact sequences
$$
\begin{array}{cccccccc}
  &  &  &  0 &    & 0  &    &  \\
  &  &  & \downarrow  &    &\downarrow  &    &  \\
       &  &  & k & \cong   & k &  &    \\
  &  &  & \downarrow  &    &\downarrow &    &  \\
      \rightarrow & \H^{0}(A, \Omega_{A}^{1}) & \longrightarrow &
\H^{1}(A, B_{i}) & \longrightarrow & \H^{1}(A, Z_{i}) & \longrightarrow &
\H^{1}(A, \Omega_{A}^{1})    \\
  & \rmapdown{=}  &  & \rmapdown{C} &    &\rmapdown{C}   &    & \rmapdown{=}   \\
      \rightarrow & \H^{0}(A, \Omega_{A}^{1}) & \longrightarrow &
\H^{1}(A, B_{i-1}) & \longrightarrow & \H^{1}(A, Z_{i-1}) & \longrightarrow &
\H^{1}(A, \Omega_{A}^{1}).   
\end{array}
$$
If $\H^{1}(A, Z_{i}) \rightarrow \H^{1}(A, Z_{i -1})$ is surjective, then
by this diagram $\H^{1}(A, B_{i}) \rightarrow \H^{1}(A, B_{i -1})$ is 
surjective. Therefore, again by this diagram we have $\H^{1}(A, B_{i})$ $=
\H^{1}(A, B_{i -1}) +1$, which contradicts Lemma \ref{lm:Bi-0}.
Therefore, 
$\H^{1}(A, Z_{i})$$\rightarrow \H^{1}(A, Z_{i -1})$ is not surjective for 
$i \geq 2$. Hence by Lemmas \ref{lm:surjective}, \ref{lm:hikaku} and 
$\dim \H^{1}(A, \Omega_{A}^{1}) = 4$, we conclude $\dim \H^{1}(A, Z_{i}) = 4$.}

\section{The de Rham Cohomology}
The de Rham cohomology of an abelian surface $A$ is the hypercohomology
 of the complex $(\Omega_A^{\bullet},d)$. On $\H^2_{dR}$ we have a perfect
pairing $\langle \, , \, \rangle$ given by Poincar\'e duality. 
The Hodge spectral sequence with $E_1^{ij}=\H^j(A, \Omega^i)$ converges
to $\H_{dR}^*(A)$. The second spectral sequence of hypercohomology has
$E_2^{ij}= \H^i({\cal H}^j(\Omega^{\bullet}))$ abutting to
$\H_{dR}^{i+j}(A/k)$. But the Cartier operator yields an isomorphism
of sheaves  
$$
C^{-1}: \Omega^i_{A^{(p)}} \stackrel{\sim}{\rightarrow} {\cal
H}^i(F_*(\Omega_{A/k}^{\bullet})),
$$
so that we can rewrite this as
$$
E_2^{ij}=\H^i(A^{\prime}, {\cal H}_j(\Omega^{\bullet}))= \H^i(A^{\prime},
\Omega_{X^{\prime}}^j) \Rightarrow \H_{dR}^*(A),
$$
where $A^{\prime}$ is the base change of $A$ under Frobenius. We thus
get two filtrations on the de Rham cohomology: the Hodge filtration
$$
(0) \subset F^2 \subset F^1 \subset \H_{dR}^2,
$$
and the conjugate filtration
$$
(0) \subset G_1 \subset G_2 \subset \H_{dR}^2.
$$
We have $\rank (F^1)=\rank (G_2)= 5$, $\rank (F^2)=\rank (G_1)=1$ and
$$
(F^1)^{\perp}=F^2\quad{\mbox{and}}\qquad G_1^{\perp}=G_2.
$$
We have also
$$
    F^1/F^2 \cong \H^{1}(A, \Omega^{1}_{A}).
$$
Note that the image of
$\H^2(A,O_A)$ under Frobenius is $G_1$. 

With respect to a suitable affine open covering $\U = \{U_{i}\}_{i\in I}$ of
$A$, we have  the following 2-cocycles which give a basis of $\H^{1}(A,
\O_{A})$ in each case.

\begin{itemize}
\item[{\rm Case}~$({\rm i})$] $A$ is ordinary.\\
$\{f_{ij}\}$, $\{g_{ij}\}$ such that
$\{f^{p}_{ij}\} \sim \{f_{ij}\}$, $\{g^{p}_{ij}\} \sim \{g_{ij}\}$. 
Therefore, there exist $f_{i}, g_{i} \in \Gamma (U_{i}, \O_{A})$ such that 
$f^{p}_{ij} = f_{ij} + f_{j} - f_{i}$ and $g^{p}_{ij} =  g_{ij} + g_{j} -
g_{i}$.
\item[{\rm Case}~$({\rm ii})$] $A$ is  of $p$-rank 1. \\
$\{f_{ij}\}$, $\{g_{ij}\}$ such that
$\{f^{p}_{ij}\} \sim 0$ and $\{g^{p}_{ij}\} \sim \{g_{ij}\}$. 
Therefore, there exist $f_{i}, g_{i} \in \Gamma (U_{i}, \O_{A})$ such that 
$f^{p}_{ij} = f_{j} - f_{i}$ and $g^{p}_{ij} = g_{ij} + g_{j} - g_{i}$ and
$\omega_{1} = df_{i}$ gives a non-zero regular 1-form on $A$.
\item[{\rm Case}~$({\rm iii})$]  $A$ is supersingular and not superspecial.\\
$\{f_{ij}\}$, $\{f^{p}_{ij}\}$ such that
$\{f^{p^{2}}_{ij}\} \sim 0$ and that $\{f^{p}_{ij}\}$ is not  cohomologous
to zero.  Therefore, there exists $f_{i} \in \Gamma (U_{i}, \O_{A})$ such
that 
$f^{p^{2}}_{ij} = f_{j} - f_{i}$. Now $\omega_{1} = df_{i}$ gives a non-zero 
regular 1-form on $A$.
\item[{\rm Case}~$({\rm iv})$] $A$ is superspecial.\\
$\{f_{ij}\}$, $\{g_{ij}\}$ such that
$\{f^{p}_{ij}\} \sim 0$, $\{g^{p}_{ij}\} \sim 0$. Therefore,
there exist $f_{i}, g_{i} \in \Gamma (U_{i}, \O_{A})$ such that 
$f^{p}_{ij} = f_{j} - f_{i}$ and $g^{p}_{ij} =  g_{j} - g_{i}$.
$\omega_{1} = df_{i}$ and $\omega_{2} = dg_{i}$ give linearly independent 
non-zero regular 1-forms on $A$. A non-zero regular 2-form on $A$ is 
given by $\omega_{1} \wedge \omega_{2}$.
\end{itemize}

By the definition of cup product for {\v C}ech cocycles, a basis  of
$\H^{2}(A, \O_{A})$ is given by $\{f_{ij}g_{jk}\}$ and the action of
Frobenius map on 
$\H^{2}(A, \O_{A})$ is given by $\{f_{ij}g_{jk}\} \mapsto 
\{f^{p}_{ij}g^{p}_{jk}\}$.  Now we consider 
the Frobenius map on the de Rham cohomology $\H^{2}_{dR}(A)$. 
On $F^{1}$ the Frobenius map $F$ is zero and so  it induces a homomorphism
$$
 F  : \H^{2}(A, \O_{A})\cong \H^{2}_{dR}(A)/F^{1} \longrightarrow 
\H^{2}_{dR}(A).
$$
The image of $F$ coincides with $G_{1}$, and 
a basis of its image is given by $(\{f^{p}_{ij}g^{p}_{jk}\}, 0, 0)$.
We denote this element by $\alpha$. The following lemma follows immediately
from $\Omega_A^1\cong O_A \oplus O_A$.

\begin{lm}\label{lm:one-one-form} Under the identification
$\H^{1}(A, \Omega^{1}_{A}) \cong \H^{1}(A, \O_{A}) \otimes_{k} 
\H^{0}(A, \Omega^{1}_{A})$
 a basis $\langle \omega_{1},  \omega_{2} \rangle$ of $\H^{0}(A,
\Omega^{1}_{A})$, defines a basis of $\H^{1}(A, \Omega^{1}_{A})$  by 
$$
    \langle  \{f_{ij}\omega_{1}\},  \{g_{jk}\omega_{1}\}, 
\{f_{ij}\omega_{2}\},  \{g_{jk}\omega_{2}\}\rangle.
$$
\end{lm}

The following result can be found in  \cite{Og1}. We give here an 
elementary proof.

\begin{lm}
An abelian surface $A$ over $k$ is superspecial if and only if the  Hodge
filtration 
$F^{2} = G_{1}$.
\end{lm}
\proof{If $A$ is ordinary, i.e., in Case $({\rm i})$, the Frobenius map 
on $\H^{2}(A, \O_{A})$ is bijective.
Therefore, $G_{1}$ is not contained in $F^{1}$. Since $F^{2} \subset F^{1}$,
we have $F^{2} \neq G_{1}$. In Case $({\rm ii})$, we have
$$
\begin{array}{cl}
  f^{p}_{ij}g^{p}_{jk} & =  (f_{j} - f_{i})g^{p}_{jk} \\
                       & = f_{j}g^{p}_{jk} - f_{i}(g^{p}_{ik} - g^{p}_{ij}).
\end{array}
$$
Therefore, considering the  element $\beta = (\{f_{i}g^{p}_{ij}\}, 0)$,
we have 
\begin{equation}
   \alpha = - \delta (\beta) +  (0, \{g^{p}_{ij}\omega_{1}\}, 0).
\label{eqn:1}
\end{equation}
Here, $\delta$ means the differential in the sense of {\v C}ech cocycles.
Hence, $\alpha = (0, \{g^{p}_{ij}\omega_{1}\}, 0)$ in $\H^{2}_{dR}(A)$.
If $(0, \{g^{p}_{ij}\omega_{1}\}, 0) \sim (0, 0, \{\Omega_{i}\}) \in F^{2}$, 
then there exists an element $(0, 0, \{\eta_{i}\})$ with 
$\eta_{i} \in \Gamma(U_{i}, \Omega^{1}_{A})$  such that
$$
 (0, \{g^{p}_{ij}\omega_{1}\}, 0) =  (0, 0, \{\Omega_{i}\})  +
\delta ((0, \{\eta_{i}\})).
$$
Therefore, we have
$$
    g^{p}_{ij}\omega_{1} = \eta_{j} - \eta_{i},~ \Omega_{i} = d\eta_{i}.
$$
By $g^{p}_{ij} = g_{ij} + g_{j} - g_{i}$ we have
$$
  g_{ij}\omega_{1} = (\eta_{j} - g_{j}\omega_{1}) - 
(\eta_{i} - g_{i}\omega_{1}).
$$ 
Therefore, $\{ g_{ij}\omega_{1}\}$ is zero in $\H^{1}(A, \Omega^{1}_{A})$, 
which contradicts Lemma \ref{lm:one-one-form}. Therefore, $\alpha$ is not
contained in $F^{2}$. Hence, we have $F^{2} \neq G_{1}$ in Case (ii).
We can prove $F^{2} \neq G_{1}$ for Case $({\rm iii})$ by the same way as 
in Case $({\rm ii})$.
Now, we consider Case $({\rm iv})$. The calculation is completely parallel 
to Case $({\rm ii})$ until the equation (\ref{eqn:1}). We take the element 
$\gamma = (0, \{g_{i}\omega_{1}\})$. Then, we have
$$
 (0, \{g^{p}_{ij}\omega_{1}\}, 0) = \delta (\gamma) + 
(0, 0 , \omega_{1}\wedge \omega_{2})
$$
by $g^{p}_{ij} =  g_{j} - g_{i}$. Therefore, we have 
$\alpha \sim (0, 0 , \omega_{1}\wedge \omega_{2}) \in F^{2}$.
Hence, we have $F^{2} = G_{1}$ in this case. This completes the proof.}

\section{Subspaces of $\H^{1}(A, \Omega_{A}^{1})$}
For an abelian surface $A$, we consider the natural inclusions
$$
B_i \hookrightarrow \Omega_{A}^{1},~  Z_i \hookrightarrow
\Omega_{A}^{1}.
$$
These inclusions induce homomorphisms
$$
\H^1( B_i) \rightarrow \H^1(\Omega_{A}^{1}),~  
\H^1(Z_i) \rightarrow \H^1(\Omega_{A}^{1}).
$$
We denote by ${\Im} \H^1(B_i)$ and ${\Im} \H^1(Z_i)$,  respectively,
the images of these homomorphisms.
In this section, we determine the image $\Im \H^{1}(A, B_{i})$ 
(resp. $\Im \H^{1}(A, Z_{i})$) of $\H^{1}(A, B_{i})$ (resp. $\H^{1}(A, Z_{i})$)
in $\H^{1}(A, \Omega_{A}^{1})$. For the proof of the following lemma, 
see \cite{GK}, Lemma 9.3.

\begin{lm}\label{lm:orthogonal}
${\Im} \H^1(B_i)$ and ${\Im} \H^1(Z_i)$ are orthogonal subspaces.  
\end{lm}

\begin{lm}
Assume $h = 1$. Then we have
$\Im \H^{1}(A, B_{i}) = 0$ and moreover
$\Im \H^{1}(A, Z_{i}) \cong \H^{1}(A, \Omega_{A}^{1})$ for $ i \geq 1$.
\end{lm}
\proof{This follows from $\H^{\ell}(A,  d\Omega_{A}^{1}) = 0$ for any $\ell$, and
the exact sequence
$$
  0 \rightarrow Z_{i+1} \longrightarrow Z_{i} 
\stackrel{dC^{i}}{\longrightarrow}
 d\Omega_{A}^{1}  \rightarrow 0.
$$
}

\begin{lm}
Assume $h = 2$. Then the subspace 
$\Im \H^{1}(A, B_{i}) \subset \H^{1}(A, \Omega_{A}^{1})$ has dimension
$1$ for any $ i \geq 1$. Moreover,
we have $\Im \H^{1}(A, B_{i}) = k\{g_{ij}\omega_{1}\}$, using the notation 
in Section 5.
\end{lm}
\proof{We consider natural inclusions 
$B_{i} \rightarrow B_{i +1} \rightarrow \Omega_{A}^{1}$.
By the exact sequence
$$
    0 \rightarrow B_{i} \longrightarrow B_{i + 1} 
\stackrel{C^{i}}{\longrightarrow} B_{1} \rightarrow 0,
$$
we get an exact sequence
$$
  0 \rightarrow \H^{0}(A, B_{1}) \longrightarrow   \H^{1}(A, B_{i}) 
\stackrel{\varphi_{i}}{\longrightarrow}  \H^{1}(A, B_{i + 1}) 
\longrightarrow 
\H^{1}(A, B_{1}).
$$
Therefore, by Lemma \ref{lm:Bi-0}, we see
that $\dim \Im \varphi_{i} = 1$. Hence, we have
$\dim \Im [\H^{1}(B_{i})  \rightarrow \H^{1}(\Omega_{A}^{1})] \leq 1$.
On the other hand, we consider the natural homomorphism
$\H^{1}(A, B_{1}) \rightarrow \H^{1}(A, \Omega_{A}^{1})$. Using the notation
in Section 5, we see, by Lemma \ref{lm:one-one-form} and 
$\{g^{p}_{ij}\} \sim \{g_{ij}\}$, that $g^{p}_{ij}\omega_{1} = d(g^{p}_{ij}f_{i})$
is a non-zero element of $\H^{1}(A, B_{1})$ whose image is not zero in 
$\H^{1}(A, \Omega_{A}^{1})$. Using natural homomorphisms
$$
  \H^{1}(A, B_{1}) \rightarrow \H^{1}(A, B_{i}) \rightarrow \H^{1}(A, 
\Omega_{A}^{1}),
$$
we see  $\dim \Im [\H^{1}(B_{i})  \rightarrow \H^{1}(\Omega_{A}^{1})]  \geq
1$. Hence, we have $\dim \Im [\H^{1}(B_{i})  \rightarrow
\H^{1}(\Omega_{A}^{1})] = 1$.}

We now consider the homomorphisms
$$
\begin{array}{cccccc}
 c_{1}^{(i)} : &\Pic (A) & = & \H^{1}(A, \O^{*}_{A}) & \longrightarrow & 
\H^{1}(A, Z_{i}) \\
                  &  D &=         & \{h_{ij}\}     &   \mapsto      &   \{
d\log {h_{ij}}\} .
\end{array}
$$
If a divisor $D$ is algebraically equivalent to a divisor $E$, then 
there  exists an element $a$ of $A$ such that $T_{a}^{*}E \sim
D$ (linearly equivalent) with $T_a$ the translation by $a$. Since $A$ is a
complete variety, the linear  representation of $A$ is trivial. Therefore,
$A$ acts on $\H^{1}(A, Z_{i})$ trivially. Therefore, we have $c_{1}^{(i)}(D)
= c_{1}^{(i)}(T_{a}^{*}E) = c_{1}^{(i)}(E)$. Hence, we have the induced
homomorphisms
$$
\begin{array}{cccc}
     c_{1}^{(i)} :  & NS(A)& \longrightarrow & \H^{1}(A, Z_{i})  \\
             & D =  \{h_{ij}\}     &   \mapsto      &   \{d \log {h_{ij}}\} .  
\end{array}
$$
Note that $c_{1} = c_{1}^{(0)}  :  NS(A)  \longrightarrow  \H^{1}(A, 
\Omega_{A}^{1})$ gives the usual Chern class.
We denote by $c_{1}^{dR}$ the Chern mapping $NS(A) \longrightarrow
H^{2}_{dR}(A)$.

The following lemma is due to Ogus \cite{Og1}.
\begin{lm} The abelian surface
$A$ is superspecial if and only if we have
$\Im c_{1}^{dR} \cap F^{2} \neq \{0\}$.
\end{lm}

\begin{lm}\label{lm:c1injection} 
Assume that $A$ is not a superspecial $($resp. superspecial$)$ abelian 
surface. Then the homomorphism  
$c_{1}^{(n)} :   NS(A)/pNS(A) \longrightarrow  \H^{1}(A, Z_{n})$
is injective for any $n \geq 0$ (resp.\ $n \geq 1$).
\end{lm}
\proof{Assume that $A$ is not superspecial.
It suffices to prove 
$$
c_{1} :  NS(A)/pNS(A)  \longrightarrow  \H^{1}(A, \Omega_{A}^{1})
$$
is injective.
Suppose there exists an element $D = \{h_{ij}\} \in NS(A)/pNS(A)$ such that
$\{d\log h_{ij}\}$ is zero in $\H^{1}(A, \Omega_{A}^{1})$. Then, there exists
an element $\omega_{i} \in \Gamma(U_{i}, \Omega_{A}^{1})$ such that
$d\log h_{ij}= \omega_{j} - \omega_{i}$. If $d \omega_{i} \neq 0$, then
we have $\Im c_{1}^{dR} \cap F^{2} \neq \{0\}$, which contradicts our 
assumption. Therefore, we have  $d \omega_{i} = 0$. Applying the Cartier
operator to both sides, we have 
$d\log h_{ij}= C(\omega_{j}) - C(\omega_{i})$. Therefore, we have 
$C(\omega_{j}) -\omega_{j} = C(\omega_{i})  - \omega_{i}$. Since
$C - {\rm id} : \H^{0}(A, \Omega_{A}^{1}) \rightarrow \H^{0}(A, 
\Omega_{A}^{1})$ is surjective, there exists an element $\omega \in
\H^{0}(A, \Omega_{A}^{1})$ such that $C(\omega_{i})  - \omega_{i} =
C{\omega} - \omega$. Replacing $\omega_{i}$ by $\omega_{i}  - \omega$, we
have
$d\log h_{ij}= \omega_{j} - \omega_{i}$ with $C(\omega_{i}) = \omega_{i}$.
Therefore, there exists $h_{i} \in \Gamma(U_{i}, \O_{A}^{*})$ such that
$\omega_{i} = d\log h_{i}$ and we see $D \in pNS(A)$.

If $A$ is superspecial it suffices to prove that
$c_{1}^{(1)} :  NS(A)/pNS(A)  \longrightarrow  \H^{1}(A, Z_{1})$ is 
injective.  The proof is completely similar to the later part of the above
proof.}

\begin{lm}\label{lm;Zinjective} 
The natural homomorphism 
$\H^{1}(A, Z_{1}) \longrightarrow \H^{2}_{dR}(A)$ is injective.
\end{lm}
\proof{Suppose there exists a {\v C}ech cocycle 
$\alpha = \{\omega_{ij}\} \in \H^{1}(A, Z_{1})$
with respect to an affine open covering $\{U_{i}\}$ of $A$ 
such that  $\alpha$ is zero in $\H^{2}_{dR}(A)$. Then by the definition
of de Rham cohomology, there exists
$\omega_{i} \in \Gamma (U_{i}, \Omega_{A}^{1})$ such that 
$\omega_{ij} = \omega_{j} - \omega_{i}$ on $U_{i} \cap U_{j}$ and 
that $d\omega_{i} = 0$. This means that $\alpha = 0$ in $\H^{1}(A, Z_{1})$.} 

The following proposition follows from Lemmas \ref{lm:c1injection} and 
\ref{lm;Zinjective}.

\begin{prop}\label{prop:Chern}
Let $A$ be an abelian surface over $k$. Then the Chern mapping
$$
      c_{1}^{dR} : NS(A)/pNS(A) \longrightarrow \H^{2}_{dR}(A)
$$
is injective.
\end{prop}
\section{The Supersingular Case}
If $A$ is a supersingular abelian surface we denote by $\sigma_{0}$ 
the Artin-invariant of $A$. This is defined as $\sigma_0= {\rm
ord}_p(\Delta)$, where $\Delta$ is the discriminant of the intersection form
$\langle \, ,\, \rangle$ on $NS(A)$. Then, $\sigma_{0}$ is equal to either 1
or 2. Moreover, $A$ is superspecial if and  only if $\sigma_{0} = 1$. We
denote by $\langle \Im c_{1} \rangle$ the subspace of $\H^{1}(A,
\Omega_{A}^{1})$ generated by $\Im c_{1}$ over $k$.
\begin{lm}\label{lm:sigma}
Assume that $A$ is a superspecial abelian surface. Then 
$$
 \langle \Im c_{1} \rangle = \H^{1}(A, \Omega_{A}^{1}).
$$
\end{lm}
\proof{For a superspecial abelian surface the endomorphism ring is
$M_2(R)$ with $R$ an order in a quaternion division algebra over ${\bf Q}$.
It follows that the rank of $NS(A)$ equals $6$. From $\sigma_0=1$ it follows
that the image of $c_1^{dR}$ has rank at least $5$ in $H^2_{dR}(A)$. Since 
$\langle \Im c_{1}^{dR} \rangle \subset F^{1}$ and $\dim F^{1} = 5$, we have 
$\langle \Im c_{1}^{dR}\rangle = F^{1}$.  Hence, the result follows from
$F^{1}/F^{2} \cong \H^{1}(A, \Omega_{A}^{1})$.}

\begin{cor}\label{cor:superspecialZi}
Assume $A$ is superspecial. Then the natural mapping 
$$
\H^{1}(A, Z_{i}) \rightarrow \H^{1}(A, Z_{i-1})
$$ 
is surjective for $i \geq 1$.
In particular, for $i \geq 1$ we have the equality $\Im \H^{1}(A, Z_{i}) =
\H^{1}(A,
\Omega_{A}^{1})$.
\end{cor}
\proof{By Lemma \ref{lm:surjective}, it suffices to show this corollary
for $\H^{1}(A, Z_{1}) \rightarrow \H^{1}(A, \Omega_{A}^{1})$.
We have 
$\langle \Im c_{1} \rangle = \H^{1}(A, \Omega_{A}^{1})$ by the assumption,
and $\Im c_{1} \subset \Im \H^{1}(A, Z_{1})$. The result follows from these facts.}

\begin{cor}
Assume $A$ is superspecial. Then, $\dim \Im \H^{1}(A, B_{i}) = 0$  for any
$i \geq 1$. Moreover, $\dim \H^{1}(A, Z_{i}) = 4 + i$.
\end{cor}
\proof{The first statement follows from Corollary
\ref{cor:superspecialZi}  and Lemma \ref{lm:orthogonal}.
For the second, note that the exact sequence
$$
      0 \rightarrow d\O_{A} \longrightarrow Z_{i+1} \stackrel{C}{\longrightarrow}Z_{i} 
 \rightarrow 0
$$
leads to the exact sequence
$$
0 \rightarrow k \longrightarrow \H^{1}(A, Z_{i+1})
\stackrel{C}{\longrightarrow}
 \H^{1}(A, Z_{i}) \longrightarrow  k.
$$
Since we have $\langle \Im c_{1} \rangle = \H^{1}(A, \Omega_{A}^{1})$, we 
have $\langle C(c^{(1)}(NS(A)))\rangle  = \H^{1}(A, \Omega_{A}^{1})$ by Lemma
\ref{lm:sigma}. Therefore,  $\H^{1}(A, Z_{1}) \longrightarrow \H^{1}(A,
\Omega_{A}^{1})$  is surjective.
Hence, by Lemma \ref{lm:surjective}, $\H^{1}(A, Z_{i+1}) 
\longrightarrow \H^{1}(A, Z_{i})$
 is surjective for any $i \geq 0$, and
we have the equality $\dim \H^{1}(A, Z_{i+1})  = \dim \H^{1}(A, Z_{i})  + 1$. 
The result follows from this fact by induction.}

\begin{lm}\label{lm:Z1}
Assume that $A$ is a supersingular abelian surface which is not superspecial.
Then $\dim \Im \H^{1}(A, Z_{1}) \geq 3$.
\end{lm}
\proof{We consider the exact sequence
$$
 0 \rightarrow d\O_{A} \longrightarrow Z_{1} 
\stackrel{C}{\longrightarrow} \Omega_{A}^{1} \rightarrow 0.
$$
By this exact sequence, we have an exact sequence
$$
 0 \rightarrow k \longrightarrow \H^{1}(A, Z_{1}) \stackrel{C}{\longrightarrow} 
\H^{1}(A,\Omega_{A}^{1}) \longrightarrow k.
$$
Since $\H^{1}(A,\Omega_{A}^{1}) = 4$, we have $\dim \Im \H^{1}(A, Z_{1}) \geq 3$.}

\begin{lm}\label{lm:B1dim}
If $A$ is a supersingular abelian surface which is not superspecial then
$\dim \Im \H^{1}(A, B_{1}) = 1$. Moreover, we have 
$\Im \H^{1}(A, B_{1}) = k\{f_{ij}^{p}\omega_{1}\}$, using the notation 
in Section 5.
\end{lm}
\proof{By Lemmas \ref{lm:Z1} and \ref{lm:orthogonal}, we have 
$\dim \Im \H^{1}(A, B_{1}) \leq 1$. Using the notation in Section 5,  
we have a non-zero element 
$f_{ij}^{p}\omega_{1} =  d(f_{ij}^{p}f_{i})$ in $\H^{1}(A, B_{1})$ which is
not zero in  $\H^{1}(A,\Omega_{A}^{1})$ by Lemma \ref{lm:one-one-form}.
Therefore, we have $\dim \Im \H^{1}(A, B_{1}) \geq 1$.
Hence, we have $\dim \Im \H^{1}(A, B_{1}) = 1$.}

\begin{cor}\label{cor:superZ}
If $A$ is a supersingular abelian surface which is not superspecial
then $\dim \H^{1}(A, Z_{1}) = 4$, and $\dim \Im \H^{1}(A, Z_{1}) = 3$.
\end{cor}
\proof{The latter part follows from Lemma \ref{lm:Z1}, 
$\dim \Im \H^{1}(A, B_{1}) = 1$
and $\Im \H^{1}(A, Z_{1})  \subset \Im \H^{1}(A, B_{1})^{\perp}$. The former part
follows from the latter part and the exact sequence
$$
   0 \rightarrow Z_{1} \longrightarrow \Omega_{X}^{1} \stackrel{d}{\longrightarrow} 
d\Omega_{X}^{1} \rightarrow 0.
$$}

Let $A$ be a supersingular abelian surface, and
let $\ell_{n}$ be the smallest non-negative integer $\ell$ such that there
exist elements 
$D_{1}, \ldots, D_{\ell}$ in  $NS(A)/pNS(A)$ with $c_{1}^{(n)}(D_{1}),
\ldots, c_{1}^{(n)}(D_{\ell})$  linearly independent over ${\bf F}_{p}$
and linearly dependent over k in $\H^{1}(A, Z_{n})$. Since we have, by Lemma 
\ref{lm:c1injection},  
$$
\dim_{{\bf F}_{p}}c_{1}^{(1)}(NS(A)/pNS(A)) = 6~ \mbox{and}~ 
\dim_{k}\H^{1}(A, Z_{1}) \leq 5,
$$
we see that $\ell_{1} \geq 1$.

\begin{lm}\label{lm:small}
Let $A$ be a supersingular abelian surface. For $n \geq 1$ we have
$\ell_{n + 1} > \ell_{n}$ if $\ell_{n + 1} \neq 0$.
\end{lm}
\proof{Assume $\ell_{n + 1} \neq 0$.
Then, there exist elements $D_{1}, \ldots, D_{\ell_{n+1}}$
$\in$ $NS(A)/pNS(A)$ such that $c_{1}^{(n+1)}(D_{1}), \ldots, c_{1}^{(n+1)}(D_{\ell_{n+1}})$
are linearly independent over ${\bf F}_{p}$, and such that 
\begin{equation}
      a_{1}c_{1}^{(n+1)}(D_{1}) +  \ldots 
+ a_{\ell_{n+1}} c_{1}^{(n+1)}(D_{\ell_{n+1}}) = 0 ~\mbox{in}~\H^{1}(A, Z_{n+1})
 \label{eqn:b1}
\end{equation}
is a non-trivial linear relation over $k$. We may assume $a_{1} = 1$. 
By the minimality of  $\ell$  we see $a_{i}/a_{j} \notin {\bf
F}_{p}$ for any $i, j$, $i \neq j$. Applying the Cartier operater $C$ to
both sides in (\ref{eqn:b1}), we get
\begin{equation}
   a_{1}c_{1}^{(n)}(D_{1}) + a_{2}^{p^{-1}}c_{1}^{(n)}(D_{2}) + \ldots 
+ a_{\ell_{n+1}}^{p^{-1}} c_{1}^{(n)}(D_{\ell_{n+1}}) = 0 ~\mbox{in}~\H^{1}(A, Z_{n})
 \label{eqn:b2}
\end{equation}
by using  $C(c_{1}^{(n+1)}(D_{j})) = c_{1}^{(n)}(D_{j})$. Subtracting
(\ref{eqn:b2}) from (\ref{eqn:b1}), we get
a non-trivial linear relation over $k$ whose length is smaller than or equal to
${\ell_{n+1}} - 1$.
Hence we have $\ell_{n + 1} > \ell_{n}$.}

\begin{cor}
For a supersingular abelian surface $A$,
 there exists a integer $n \geq 2$ such that the natural
homomorphism
\begin{equation}
   NS(A)/pNS(A) \otimes_{{\bf F}_{p}} k \rightarrow \H^{1}(A, Z_{n})
 \label{eqn:b3} 
\end{equation}
is injective. 
\end{cor}
\proof{Since $\dim_{{\bf F}_{p}}NS(A)/pNS(A) = 6$, we have
$\ell_{n} \leq 6$. Hence, by Lemma \ref{lm:small}, there exists a positive 
integer $n$ such that  $\ell_{n} = 0$.}

\begin{lm}\label{lm:contain}
Assume that $A$ is a supersingular abelian surface which is not superspecial.
Let $s$ be the smallest integer such that $($\ref{eqn:b3}$)$ is injective
and let $\varphi_{t}^{s} : \H^{1}(A, B_{t}) \rightarrow \H^{1}(A, Z_{s})$ be
the homomorphism  induced by the natural inclusion 
$B_{t} \rightarrow Z_{s}$. Then 
$\Im \varphi_{t}^{s} \subset \langle c_{1}^{(s)}(NS(A)/pNS(A))\rangle$ for any
$t \geq 1$.
\end{lm}
\proof{Since $\dim_{{\bf F}_{p}}NS(A)/pNS(A) = 6$ and 
$\dim \H^{1}(A, Z_{1}) = 4$ by Corollary \ref{cor:superZ}, 
we see $s \geq 3$. 
We consider an exact sequence
$$
   0 \rightarrow B_{1} \longrightarrow Z_{s} \stackrel{C}{\longrightarrow} Z_{s-1}
\rightarrow 0.
$$
By the long exact sequence, we have $\dim \Im \varphi_{1}^{s} = 1$.

By the assumption on $s$ 
there exist elements $D_{1}, \ldots, D_{\ell_{s-1}}$
in the ${\bf F}_p$-vector space $NS(A)/pNS(A)$ such that 
$c_{1}^{(s-1)}(D_{1}), \ldots, c_{1}^{(s-1)}(D_{\ell_{s-1}})$
are linearly independent over ${\bf F}_{p}$, and such that 
\begin{equation}
     a_{1}c_{1}^{(s-1)}(D_{1}) +  \ldots 
+ a_{\ell_{s-1}} c_{1}^{(s-1)}(D_{\ell_{s-1}}) = 0 ~\mbox{in}~\H^{1}(A, Z_{s-1})
 \label{eqn:b4}
\end{equation}
is a non-trivial linear relation over $k$. We may assume $a_{1} = 1$. 
Then, by the smallestness of $\ell_{s-1}$, we see 
$a_{i}/a_{j} \notin {\bf F}_{p}$ for any $i, j$, $i \neq j$.
We apply the Cartier inverse to (\ref{eqn:b4}) and find
$$
 a_{1}c_{1}^{(s)}(D_{1}) + a_{2}^{p}c_{1}^{(s)}(D_{2}) + \ldots 
+ a_{\ell_{s-1}}^{p} c_{1}^{(s)}(D_{\ell_{s-1}}) = \{dh_{ij}\} ~{\mbox{in}}~\H^{1}(A, Z_{s})
$$
with $\{dh_{ij}\} \in \H^{1}(A, B_{1})$. 
If $\{dh_{ij}\} = 0$ in $\H^{1}(A, B_{1})$
we see that 
$$
NS(A)/pNS(A) \otimes_{{\bf F}_{p}} k \rightarrow \H^{1}(A, Z_{s})
$$
is not injective, a contradiction. Therefore,  $\{dh_{ij}\} \neq 0$ in 
$\H^{1}(A, B_{1})$, and $\{dh_{ij}\}$ gives a basis of $\Im \varphi_{1}^{s}$.
So we have $\Im \varphi_{1}^{s} \subset \langle
c_{1}^{(s)}(NS(A)/pNS(A))\rangle$.

Suppose that there exists a $j$ such that $\Im \varphi_{j}^{s}$ is not contained
in k-vector space $\langle c_{1}^{(s)}(NS(A)/pNS(A))\rangle$. We take 
the smallest such $j$.
Then, we have $j \geq 2$ as we showed above. We take an element 
$\alpha \in \Im \varphi_{j}^{s}$ such that 
$\alpha \notin \langle c_{1}^{(s)}(NS(A)/pNS(A))\rangle$. 
Then, there exists an element $\tilde{\alpha} \in \H^{1}(A, B_{j})$
such that $\varphi_{j}^{s} (\tilde{\alpha}) = \alpha$. We have 
$C(\tilde{\alpha}) \in \H^{1}(A, B_{j-1})$ and 
$\varphi_{j-1}^{s}(C(\tilde{\alpha})) \in \langle c_{1}^{(s)}(NS(A)/pNS(A))\rangle$
by the assumption on $j$. Applying the inverse Cartier, we see that there
exists an element $\{dk_{ij}\} \in \H^{1}(A, B_{1})$ such that
$\alpha + \{dk_{ij}\}$ is contained in  $\langle c_{1}^{(s)}(NS(A)/pNS(A))\rangle$.
Since 
$$
\{dk_{ij}\} \in \H^{1}(A, B_{1}) \subset \langle c_{1}^{(s)}(NS(A)/pNS(A))\rangle,
$$
we have $\alpha \in \langle c_{1}^{(s)}(NS(A)/pNS(A))\rangle$, a contradiction.
This completes the proof.}

\begin{prop}\label{prop:Zidim}
For a supersingular abelian surface $A$ which is not superspecial
 $\H^{1}(A, Z_{i+1}) \rightarrow \H^{1}(A, Z_{i})$ is surjective for 
$i \geq 1$, and $\H^{1}(A, Z_{1}) \rightarrow \H^{1}(A, \Omega_{A}^{1})$
is not surjective. Moreover, we have  $\dim \H^{1}(A, Z_{i}) = i + 3$ and 
$\dim \Im \H^{1}(A, Z_{i}) = 3$ for $i \geq 1$. 
\end{prop}
\proof{We take the smallest integer $s$ such that (\ref{eqn:b3}) is injective.
Since the images of $NS(A)/pNS(A) \otimes_{{\bf F}_{q}} k$
by $C^{s}$ and by the natural homomorphism from $\H^{1}(A, Z_{s})$
to $\H^{1}(A, \Omega_{A}^{1})$ coincide with each other, we have 
an exact sequence
$$
\begin{array}{c}
  0 \rightarrow k^{\oplus 2} \longrightarrow  \H^{1}(A, B_{s})  \\
\longrightarrow \langle c_{1}^{(s)}(NS(A)/pNS(A))\rangle
 \stackrel{C^{(s)}}{\longrightarrow}  \langle c_{1}(NS(A)/pNS(A))\rangle \rightarrow 0.
\end{array}
$$
Here we use the result in Lemma \ref{lm:contain}.
Since we have 
$$
\begin{array}{l}
\dim \H^{1}(A, B_{s})) = s + 2, \\
\dim \langle c_{1}(NS(A)/pNS(A))\rangle = 3, \\
\dim \langle c_{1}^{(s)}(NS(A)/pNS(A)) = 6,
\end{array}
$$
we have $2 + 6 = (s + 2) + 3$.
Therefore, we have $s = 3$ and we have $\dim \H^{1}(A, Z_{3}) \geq 6$.
Now, by Lemmas \ref{lm:surjective}, \ref{lm:hikaku} and Corollary
\ref{cor:superZ},  we know that $\dim \H^{1}(A, Z_{3}) = 6$
and $\langle c_{1}^{(3)}(NS(A)/pNS(A)) \rangle = \H^{1}(A, Z_{3})$. The rest follows from 
Lemma \ref{lm:surjective}.}

\begin{cor}\label{cor:Bidim}
Assume that $A$ is a supersingular abelian surface which is not superspecial.
Then $\dim \Im \H^{1}(A, B_{i}) = 1$ for any $i \geq 1$.
\end{cor}
\proof{This follows from Lemma \ref{lm:B1dim} 
and Proposition \ref{prop:Zidim}.} 

\begin{cor}
Assume that $A$ is a supersingular abelian surface.
Then, we have $\langle c_{1}^{(i)}(NS(A)/pNS(A)) \rangle = \H^{1}(A, Z_{i})$ for
$ i = 1, 2, 3$. Moreover, $\langle c_{1}(NS(A)/pNS(A)) \rangle = \Im \H^{1}(A, Z_{i})$
for any $i$.
\end{cor}
\proof{Using Lemma \ref{lm:sigma} and Corollary \ref{cor:superspecialZi}, 
we can easily prove the case of superspecial abelian surface. 
Assume $A$ is supersingular
and not superspecial. The proof of the former part is given in the proof of 
Proposition \ref{prop:Zidim} by Lemma \ref{lm:surjective}.
Since $\langle c_{1}(NS(A)/pNS(A)) \rangle \subset \Im \H^{1}(A, Z_{i})$,
we can prove this corollary by Lemmas \ref{lm:surjective},
\ref{lm:hikaku} and Proposition \ref{prop:Zidim}.}

We summarize the results on the dimension of $\H^j(Z_i)$ in the following
table. The results about $h^0(Z_i)$ and $h^2(Z_i)$ follow from
\ref{cor:superspecialZi}, \ref{prop:Zidim} and the exact
sequence $0\to Z_i \to Z_{i-1} \to d\Omega_A^1\to 0$.

\bigskip
\font\tablefont=cmr8
\def\quad{\hskip 0.6em\relax}
\def\quod{\hskip 0.6em\relax}
\def\vhop{
    height2pt&\omit&&\omit&&\omit&&\omit&\cr}
\noindent{\bf Table.}
$$
\vcenter{
\tablefont
\lineskip=1pt
\baselineskip=10pt
\lineskiplimit=0pt
\setbox\strutbox=\hbox{\vrule height .7\baselineskip
                                depth .3\baselineskip width0pt}%
\offinterlineskip
\hrule
\halign{&\vrule#&\strut\quod\hfil#\quad\cr
\vhop
&${\rm type}$&&$h^0(Z_i)$&&$h^1(Z_i)$&&$h^2(Z_i)$&\cr
\vhop
\noalign{\hrule}
\vhop
&$h=1,2$&&$2$&&$4$&&$2$&\cr
\vhop
\noalign{\hrule}
\vhop
&$h=\infty$, $a=1$&&$2$&&$3+i$&&$1+i$&\cr
\vhop
\noalign{\hrule}
\vhop
&$h=\infty$, $a=2$&&$2$&&$4+i$&&$2+i$&\cr
\vhop
\noalign{\hrule}
\vhop
}
\hrule
}
$$

\section{
Degenerate Abelian Surfaces}
\smallskip
\noindent
The moduli space ${\cal A}_2$ can be compactified in a canonical way to
a compactification $\tilde{\cal A}_2$. This is a toroidal
compactification, called the Delaunay-Voronoi compactification. It is
a blow-up of the Satake-compactification due to Igusa, cf.\ also \cite{FC}. 
Alexeev gives a functorial description of it in \cite{Al}. The four strata in
codimension
$3$, $2$, $1$ and $0$ correspond to decompositions of the plane ${\bf R}^2$
by triangles, by $4$-gons, by infinite strips and by one big cell
covering the whole space. The corresponding degenerations are copies
of ${\bf P}^2$ glued, copies of ${\bf P}^1 \times {\bf P}^1$ glued, certain
non-normal compactifications of a $\G_m$-bundle over an elliptic curve or
honest abelian surfaces. We now first deal with the simplest type of
degeneration, the so-called rank $1$ degenerations, and these occur in
codimension~$1$.

A rank-$1$ degeneration of an abelian surface is a non-normal surface
obtained as follows. Start with a semi-abelian surface $G$, i.e.\ with
a $\G_m$-extension of an elliptic curve $E$
$$
1 \to \G_m \to G \to E \to 0 \eqno(*)
$$
and consider the associated ${\bf P}^1$-bundle $\pi: \tilde{G} \to E$. Then
$\tilde{G}-G$ is a union of two sections $\tilde{G}_0\sqcup
\tilde{G}_{\infty}$ of $\tilde{G}$ over $E$. The extension class of (*)
is an element $e \in \hat{E}$, the dual elliptic curve. The
compactification $X$ of $G$ is obtained by glueing the zero section
$\tilde{G}_0$ with the infinity section $\tilde{G}_{\infty}$ by a
translation over $e \in E \cong \hat{E}$. We then have
$\tilde{G}_0-\tilde{G}_{\infty}\equiv \pi^{-1}(p-q)$ with $p-q=e \in
E$. The divisor $\tilde{G}_0 +\pi^{-1}(0)$ defines a line bundle
$\tilde{L}$ on $\tilde{G}$ which descends to a line bundle $L$ on $X$
with $h^0(L)=1$. This defines the principal polarization on $X$, cf.\
\cite{Mu3}. The surface thus constructed are the fibres of the
universal family $\tilde{\cal X}_2$ over $\tilde{\cal A}_2$ over the
codimension $1$ stratum in the boundary.

The surface $X$ has $h^1(X,O_X)=4$, $h^2(X,O_X)=1$ and we can define a
formal Brauer group $\Phi_X$ associated to $\H^2(X,O_X)$ as we did for
abelian surfaces. This is a $1$-dimensional formal group. The main
result of this section is the following.

\begin{prop} The formal group $\Phi_X=\Phi^2_X$ of the rank $1$ degeneration
$X$ is isomorphic to the formal group $\Phi^1_E$ of the elliptic curve $E$.
\end{prop}
\proof{ One way of defining the formal Brauer group is as the functor
which associates to each nilpotent (commutative quasi-projective)
$O_X$-algebra $N$ the group
$$
\Phi_Y(N)= \H^2(Y,(1+O_X\otimes N)^{*}),
$$
cf. \cite{Sha}. We now use the Leray spectral sequence for the normalization
morphism
$f:\tilde{G} \to X$ of $X$. We have 
$$
E_2^{ij} = \H^j(X, R^if_*(F)) \qquad  {\rm with }
\quad F= (1+O_{\tilde{G}}\otimes N)^{*}.
$$
Because the fibres of $f$ are zero-dimensional $R^1f_*(F)$ and
$R^2f_*(F)$ vanish. The spectral sequence yields on level $2$:
$$
\matrix{ 0 & 0 &0 \cr 0 & 0& 0 \cr
\H^0(X,A) & \H^1(X,A) & \H^2(X,A)\cr
}
$$
with $A=R^0f_*(F)$. We have an exact sequence on sheaves on $X$
$$
0 \to F_X \longrightarrow A \longrightarrow B\to 0.
$$
Here $F_X= (1+O_X\otimes N)^*$ and the sheaf $B$ is concentrated on the
non-normal locus
$S$ of
$X$ where the two sections $\tilde{G}_0$ and $\tilde{G}_{\infty}$ are
glued. 

Let $t_e$ be the translation over $e\in E$. Then with $F_E=(1+O_E
\otimes N)^*$ there is an exact sequence
$$
0\to F_E {\buildrel \alpha \over \longrightarrow} F_E \oplus t_e^*F_E \to
B|S \to 0
$$
of sheaves on $S\cong E$. Here $\alpha(s)= (s,t_e^*s)$. This implies that
$\H^1(S,B|S)\cong \H^1(E,F_E)= \Phi^1_E(N)$. Since $\tilde{G}$ is ruled
we have $\H^2(\tilde{G}, F_{\tilde{G}})=0$, hence $\H^2(X,A)=0$. The spectral
sequence shows that $\H^1(\tilde{G},F_{\tilde{G}})\cong \H^1(X,A)$. In the
long exact sequence
$$
\to \H^1(F_X) \to \H^1(A) \to \H^1(B) \to \H^2(F_X)\to \H^2(A)=0
$$
the map $\H^1(F_X) \to \H^1(A)$  can be identified with
the map $\H^1(X,F_X) \to H^1(\tilde{G}, F_{\tilde{G}})$. Since the formal
group $\Phi^1_X$ fits into an extension 
$$
1 \to \hat{G}_m \to \Phi_X^1 \to \Phi_{\tilde{G}}^1 \to 0.
$$
we see that $\H^1(X,F_X) \to \H^1(\tilde{G},F_{\tilde{G}})$ is surjective.
This implies that
$\H^1(B)
\cong \H^2(F_X)$. }

We can do something similar for rank-2 degenerations of abelian surfaces. 
Then the formal group that we get is a multiplicative group $\hat{G}_m$.
The loci $M^{(h)}$ thus can be extended to the compactification $\tilde{\cal
A}_2$. We can also consider cohomology with coefficients in $W_i(O_X)$ for
such degenerate surfaces. The results we obtain are similar. Thus we can
extend the loci $M_F^{(h)}$ to the compactification.

The considerations on differential forms using the sheaves $\Omega^1$, $B_i$
and $Z_i$ can be extended to degenerate surfaces $X$ by replacing
differential forms by differential forms with log poles along the divisor of
non-normal points on $X$.

\section{The Classes of the Loci of Given Height}
Let ${\cal A}_{2}$ (or ${\cal A}_{2, n}$ with $n \geq 3, ~p \not| n$) be the
moduli stack (or fine moduli  space) of principally polarized abelian
surfaces (with level $n$-structure) and let
$\pi:{\cal X} \rightarrow {\cal A}_{2}$ be the universal family of
principally polarized abelian surfaces. We know that ${\cal A}_2$ 
(resp.\ ${\cal A}_{2,n}$) is a 3-dimensional algebraic stack 
(resp.\ variety) which
has
a natural compactification $\tilde{\cal A}_2$ (resp.\ ${\cal A}_{2,n}$), 
see e.g. [1]. It possesses a
universal family $\tilde{\pi} : \tilde{{\cal X}} \rightarrow \tilde{\cal A}_2$
(resp.\ $\tilde{\cal A}_{2, n}$). We set
$\tilde{M} = \tilde{\cal A}_2$ (or $\tilde{M}=\tilde{\cal A}_{2, n}$) and we
denote by
$v$ the class of  the coherent sheaf 
$R^{0}\pi_*\Omega^{2}_{\tilde{\cal X}/\tilde{M}}$ in the Chow group $CH_{\bf
Q}^{1}(\tilde{M})$. We set 
$$
M^{(h)}:=\{ s \in \tilde{M} : h(\tilde{\cal X}_s) \geq h\}.
$$ 
Then, as is well-known, $M^{(h)}$ is an algebraic subvariety (maybe reducible), 
and is of codimension $h-1$ in $\tilde{\cal A}_{2}$ for $h=1, 2$, and
of codimension 2 in $\tilde{\cal A}_{2}$ for $h= \infty$. Moreover,
$M^{(\infty)}$ is contained in ${\cal A}_{2}$ (or ${\cal A}_{2, n}$).

The direct image sheaves $R^2\pi_*W_i(\O_{\tilde{\cal X}})$ are coherent
$W_i(\O_{\tilde{M}})$ sheaves. Note that these sheaves are not 
$\O_{\tilde{M}}$-modules
if $i \geq 2$. 
We denote by ${M}_F^{(h)}$ $(h = 1, 2, 3)$ the scheme-theoretic zero locus
of the relative Frobenius map on 
$R^2\tilde{\pi}_*W_{(h-1)}(\O_{\tilde{\cal X}})$.
We set $M_F^{(\infty)}= M_F^{(3)}$. Then we have 
$({M_F^{(\infty)}})_{red} = M^{(\infty)}$.

\begin{th}
The class of $M_F^{(h)}$ for $h=2, \infty$ in the Chow group $CH_{\bf
Q}^{*}(\tilde{M})$   are given by
$$
     M_F^{(2)} = (p-1)v, ~M_F^{(\infty)} = (p-1)(p^2-1)v^{2}.
$$
Moreover, $M_F^{(2)} \backslash M_F^{(\infty)}$ is non-singular.
\end{th}
\proof{ The locus $M_F^{(2)}$ is given by the support of the cokernel
of 
$F : R^2\tilde{\pi}_{*}\O_{\tilde{\cal X}} \rightarrow 
R^2\tilde{\pi}_*\O_{\tilde{\cal X}}$. 
Since $F$ is a $p$-th linear mapping, it follows that $M_F^{(2)} = (p-1)v$.
By the same method as we used in van der Geer and 
Katsura \cite{GK} or as in Ogus \cite{Og1} we see that the tangent space of a
point  $x = (A_{0}, D_{0}) \in M_F^{(2)} \backslash M_F^{(\infty)}$
is given by $\Im \H^{1}(A_{0}, Z_{1}) \cap D_{0}^{\perp}$, which is 
2-dimensional. Hence, $M_F^{(2)} \backslash M_F^{(\infty)}$ is
non-singular.

Now, let $x = (A_{0}, D_{0})$ be a point in $M_F^{(\infty)}$ such 
that  $A_{0}$ is a non-superspecial abelian surface. 
We denote by $S$ the formal scheme around $x$ in $M_F^{(2)}$, 
and we denote by $X \rightarrow S$ the family
which is obtained by the restriction of $\tilde{\pi}:\tilde{{\cal X}}
 \rightarrow \tilde{M}$
to $S$. Then, the Frobenius mapping $F$ is zero on $\H^{2}(X, \O_{X/S})$. 
We denote by $\nabla$ the Gauss-Manin connection of $\H^{2}_{dR}(X/S)$. 
We consider the Hodge filtration 
$0 \subset F^{2} \subset F^{1} \subset \H^{2}_{dR}(X/S)$, and construct,
in the same way as in Section 8  of van der Geer and Katsura \cite{GK}, 
the homomorphism 
$$
\Phi_{2} : \H^{2}(W_{2}(O_{X/S})) \rightarrow \H^{2}_{dR}(X/S).
$$
We take a basis $\omega$ of
$\H^{0}(\Omega_{X/S}^{1})$ and take the dual basis $\zeta$ of
$\H^{2}(O_{X/S})$. We take a lifting
${\tilde \zeta} \in \H^{2}_{dR}(X/S)$ of $\zeta$. Then we
have $\langle {\tilde \zeta}, \omega\rangle = 1$. Since we have
a surjective homomorphism
$R: \H^{2}(W_{2}(O_{X/S})) \rightarrow \H^{2}(O_{X/S})$, there exists
an element $\alpha \in \H^{2}(W_{2}(O_{X/S}))$ such that $R(\alpha) = \zeta$.
We set
$$
 g_{h} = \langle \Phi_{h}(\alpha), {\omega}\rangle.
$$
Since $\Phi_{h}(\alpha) - g_{h}{\tilde \zeta}$ is orthogonal
to $\omega$, the element $\Phi_{h}(\alpha) - g_{h}{\tilde \zeta}$ is contained
in $F^{1}$.
 Therefore, using 
the natural isomorphism 
$\H^{2}_{dR}/F^{1} \cong \H^{2}(X, \O_{X})$, we conclude
that
$$
 \phi_{h}(\zeta) = g_{h}{\zeta} \quad {\mbox{in}}~\H^{2}(X, \O_{X}).
$$
This means that $g_{h} = 0$ gives the scheme
theoretic zero-locus $M_F^{(\infty)}$ of $\phi_{h}$.  By Theorem
\ref{th:h=2}, 
$({M_F^{(\infty)}})_{red}$ coincides with $M^{(\infty)}$. Since
$\phi_{h}$ is a $\sigma^{2}$-linear homomorphism, we conclude that
the class $M_F^{(\infty)} = (p^{2} - 1)v^{2}\vert_{M_F^{(2)}}$ in the
Chow ring 
$CH^{1}_{\bf Q}(M_F^{(2)})$. Hence, as in van der Geer and Katsura
\cite{GK}, using the projection formula, we complete the proof. }

By Ekedahl and van der Geer (cf. \cite{G2}), we have
$M^{(2)} = (p - 1)\lambda_1 $ and $M^{(\infty)} = (p - 1)(p^{2} -
1)\lambda_2$ in the Chow ring $CH_{\bf Q}^{*}(\tilde{M})$. Here $\lambda_i$
for $i=1,2$ are the Chern classes of the Hodge bundle
$E= R^0\pi_{*}\Omega^1_{\tilde{\cal X}/\tilde{M}}$, \cite{G1, G2}. We know the relation
$\lambda_1^2= 2\lambda_2$ in the Chow ring (\cite{G1}) and we have $v=
\lambda_1$. Thus we have the following relations between $M_F^{(h)}$
and $M^{(h)}$. Note that we know that all components of $M_F^{(\infty)}$ have
multiplicity
$\geq 2$.

\begin{th}
In the Chow ring $CH_{\bf Q}^{*}(\tilde{M})$, we have
$$
      M_F^{(2)} = M^{(2)}, ~M_F^{(\infty)} = 2M^{(\infty)}.
$$
All components of $M^{(\infty)}$ occur with  multiplicity $2$
in $M_F^{(\infty)}$.\end{th}

\vspace{0.5cm}

\begin{flushleft}
\begin{minipage}[t]{5.8cm}
G. van der Geer

Korteweg-de-Vries Instituut

Universiteit van Amsterdam

Plantage Muidergracht 24

1018 TV Amsterdam

The Netherlands

geer@wins.uva.nl 
\end{minipage}
\quad
\begin{minipage}[t]{6.0cm}
T. Katsura

Graduate School of Mathematical Sciences

University of Tokyo

3-8-1 Komaba, Meguro

Tokyo 153-8914

Japan

tkatsura@ms.u-tokyo.ac.jp
\end{minipage}
\end{flushleft}

\end{document}